\documentclass[12pt]{article}

\usepackage{graphicx,tikz,color,url}
\usepackage{amsmath,amssymb,latexsym}

\newtheorem{theorem}{Theorem}[section]

\newtheorem{lemma}[theorem]{Lemma}
\newtheorem{proposition}[theorem]{Proposition}

\newtheorem{corollary}[theorem]{Corollary}

\newtheorem{problem}[theorem]{Problem}

\newcommand{\proof}{\noindent{\bf Proof.\ }}
\newcommand{\qed}{\hfill $\square$ \bigskip}
\newcommand{\cp}{\,\square\,}
\newcommand{\cg}{\gamma_{\rm cg}}
\newcommand{\cc}{\gamma_{\rm c}}
\newcommand{\tcg}{\gamma_{\rm tcg}}
\newcommand{\ggt}{\gamma_{\rm tg}}
\newcommand{\ggts}{\gamma_{\rm tg}'}

\let\oldenumerate\enumerate
\renewcommand{\enumerate}{
  \oldenumerate
  \setlength{\itemsep}{0pt}
  \setlength{\parskip}{0pt}
  \setlength{\parsep}{0pt}
}

\begin{document}

\title{Total connected domination game}

\author{Csilla Bujt\'as$^{a}$\thanks{Email: \texttt{csilla.bujtas@fmf.uni-lj.si}}
\and Michael A. Henning $^{b}$\thanks{Email: \texttt{mahenning@uj.ac.za}}
\and Vesna Ir\v si\v c$^{a,c}$\thanks{Email: \texttt{vesna.irsic@fmf.uni-lj.si}}
\and Sandi Klav\v zar $^{a,c,d}$\thanks{Email: \texttt{sandi.klavzar@fmf.uni-lj.si}}
}

\date{\today}
\maketitle

\begin{center}
$^a$ Faculty of Mathematics and Physics, University of Ljubljana, Slovenia\\
\medskip

$^b$ Department of Mathematics and Applied Mathematics, University of Johannesburg, Auckland Park, 2006 South Africa\\
\medskip

$^c$ Institute of Mathematics, Physics and Mechanics, Ljubljana, Slovenia\\
\medskip

$^d$ Faculty of Natural Sciences and Mathematics, University of Maribor, Slovenia\\
\medskip
\end{center}

\begin{abstract}
The (total) connected domination game on a graph $G$ is played by two players, Dominator and Staller, according to the standard (total) domination game with the additional requirement that at each stage of the game the selected vertices induce a connected subgraph of $G$. If Dominator starts the game and both players play optimally, then the number of vertices selected during the game is the (total) connected game domination number ($\gamma_{\rm tcg}(G)$) $\gamma_{\rm cg}(G)$ of $G$. We show that $\gamma_{\rm tcg}(G)\in \{\gamma_{\rm cg}(G), \gamma_{\rm cg}(G) + 1, \gamma_{\rm cg}(G) + 2\}$, and consequently define $G$ as Class~$i$ if  $\gamma_{\rm tcg}(G) = \gamma_{\rm cg} + i$ for $i \in \{0,1,2\}$. A large family of Class $0$ graphs is constructed which contains all connected Cartesian product graphs and connected direct product graphs with minumum degree at least $2$. We show that no tree is Class~$2$ and characterize Class~$1$ trees. We provide an infinite family of Class~$2$ bipartite graphs.
\end{abstract}

\indent {\small \textbf{AMS Subj.\ Class.\ (2010):} 05C57, 05C69, 91A43}

%%%%%%%%%%%%%%%%%%%%%%%%%%%%%%%%%%%%%%%%%%%%%%%%%%%%%
%%%%%%%%%%%%%%%%%%%%%%%%%%%%%%%%%%%%%%%%%%%%%%%%%%%%%
\section{Introduction}
\label{sec:intro}
%%%%%%%%%%%%%%%%%%%%%%%%%%%%%%%%%%%%%%%%%%%%%%%%%%%%%
%%%%%%%%%%%%%%%%%%%%%%%%%%%%%%%%%%%%%%%%%%%%%%%%%%%%%

The domination game introduced in 2010 in~\cite{bresar-2010}, and the total domination game put forward in 2015 in~\cite{he-kl-ra-2015}, were studied in depth by now, the respective lists of papers~\cite{bujtas-2015, dorbec-2015, james-2019, kinnersley-2013, klavzar-2019, xu-2018} and~\cite{bresar-2017, bujtas-2018, dorbec-2016, henning-2017a, henning-2016, irsic-2019} form just a selection of these studies. The two games are in some respects similar, for instance, they both admit the so-called Continuation Principle, but also significantly different in other respects, say in the conjectured upper bounds on the corresponding invariants in terms of the order of a graph~\cite{henning-2017a, kinnersley-2013}.

Recently, in 2019, the connected domination game was introduced by Borowiecki, Fiedorowicz, and Sidorowicz~\cite{borowiecki-2019} and further investigated in~\cite{bujtas-2019, irsic-2019+}. Although the connected domination game mostly follows the definition of the standard domination game, the new game is significantly different. For instance, on the class of trees, the connected game domination number can be obtained in linear time, while the complexity of determining the game domination number is open and we suspect that it is NP-hard. In other cases, finding optimal stategies appears difficult in both games, but the connected one still being more accessible, as for instance on the games played on Cartesian product graphs, see~\cite{borowiecki-2019, bujtas-2019}.

In this paper we introduce the total version of the connected domination game following the above mentioned pattern of the total domination game~\cite{he-kl-ra-2015} versus the domination game~\cite{bresar-2010}. We proceed as follows. In Section~\ref{S:notation} we give additional definitions needed, in Section~\ref{sec:connected-games} the connected and total connected domination games are described and discussed. In Section~\ref{sec:relation} we prove that $\gamma_{\rm tcg}(G)\in \{\gamma_{\rm cg}(G), \gamma_{\rm cg}(G) + 1, \gamma_{\rm cg}(G) + 2\}$ and consequently partition the class of all graphs into Classes 0, 1, and 2, depending on which of the three possibilities holds. Then, in Section~\ref{sec:classes}, we determine different families of graphs belonging to respective classes. In the concluding section several open problems are listed.

\subsection{Notation}
\label{S:notation}

For graph theory notation and terminology, we generally follow~\cite{HeYe-book-2013}.  Specifically, let $G$ be a graph with vertex set $V(G)$ and edge set $E(G)$, and of order $n(G) = |V(G)|$ and size $m(G) = |E(G)|$. We denote the degree of a vertex $v$ in $G$ by $d_G(v)$. The minimum and maximum degrees among the vertices of $G$ are denoted by $\delta(G)$ and $\Delta(G)$, respectively. Two vertices are \emph{neighbors} if they are adjacent. The \emph{open neighborhood} of a vertex $v$ in a graph $G$, denoted $N_G(v)$, is the set of neighbors of $v$ in $G$, while the \emph{closed neighborhood} of $v$ is the set $N_G[v] = N_G(v) \cup \{v\}$. The \emph{open neighborhood} of a set of vertices $S \subseteq V(G)$ is $N_G(S) = \bigcup_{v \in S} N_G(v)$, and the \emph{closed neighborhood} of $S$ is $N_G[S] = \bigcup_{v \in S} N_G[v]$. If the graph $G$ is clear from context, we may omit the subscript in the above definitions.

The \emph{Cartesian product} $G\cp H$ of graphs $G$ and $H$ has the vertex set $V(G)\times V(H)$, vertices $(g,h)$ and $(g',h')$ being adjacent if either $gg'\in E(G)$ and $h=h'$, or $g=g'$ and $hh'\in E(H)$. The \emph{direct product} $G\times H$ has the same vertex set $V(G)\times V(H)$, vertices $(g,h)$ and $(g',h')$ being adjacent if $gg'\in E(G)$ and $hh'\in E(H)$~\cite{hik-2011}.

Let $G$ be a graph with $V(G) = \{v_1,\ldots, v_n\}$ and let ${\cal H} = \{H_1,\ldots, H_n\}$ be a family of pairwise vertex disjoint graphs. Then the {\em generalized corona} $G\odot {\cal H}$ is the graph obtained from the disjoint union of $G, H_1, \ldots, H_n$, by joining, for each $i\in [n]$, every vertex of $H_i$ with the vertex $v_i$ of $G$, cf.~\cite{dettlaff-2016}.  If all the graphs $H_i$ are isomorphic to a given graph $H$, then we may write $G\odot H$ instead of $G\odot {\cal H}$. In particular, $G\odot K_1$ is the {\em corona} over $G$.

A set $S \subseteq V(G)$ is a \emph{dominating set} of the graph $G$ if $N[S] = V(G)$. The minimum cardinality of a dominating set is the \emph{domination number} $\gamma(G)$ of $G$. A \emph{connected dominating set} is a dominating set  $S$ with the additional property that the subgraph $G[S]$ induced by $S$ is connected. The minimum cardinality of a connected dominating set is the \emph{connected domination number} $\cc(G)$ of $G$. Note that $\cc$ is defined only for connected graphs. 

A \emph{total dominating set} of a graph $G$ is a set $S$ of vertices of $G$ such that every vertex has a neighbor in $S$.  The \emph{total domination number} of an isolate-free graph $G$, denoted by $\gamma_t(G)$, is the minimum cardinality of a total dominating set in $G$. A vertex \emph{totally dominates} a vertex if they are neighbors, that is, the vertices totally dominated by a vertex $v$ are the vertices that belong to the open neighborhood $N(v)$ of $v$. We note that if $S$ is a total dominating set of $G$, then every vertex in $G$ is totally dominated by at least one vertex in~$S$.

Finally, for an integer $k \ge 1$ we denote $[k] = \{1,\ldots,k\}$.

%%%%%%%%%%%%%%%%%%%%%%%%%%%%%%%%%%%%%%%%%%%%%%%%%%%%%
%%%%%%%%%%%%%%%%%%%%%%%%%%%%%%%%%%%%%%%%%%%%%%%%%%%%%
\subsection{Connected domination games}
\label{sec:connected-games}
%%%%%%%%%%%%%%%%%%%%%%%%%%%%%%%%%%%%%%%%%%%%%%%%%%%%%
%%%%%%%%%%%%%%%%%%%%%%%%%%%%%%%%%%%%%%%%%%%%%%%%%%%%%

In this section we recall the connected domination game and the connected domination game with Chooser, and introduce the total connected domination game.

The {\em connected domination game} on a (connected) graph $G$ is played by Dominator and Staller. They play in turns, at each move selecting a single vertex of $G$ such that it dominates at least one vertex that is not yet dominated with the previously played vertices and such that at each stage of the game the selected vertices induce a connected subgraph of $G$. If Dominator has the first move, then we speak of a {\em D-game}, otherwise they play an {\em S-game}. When the game is finished, that is, when there is no legal move available, the players have determined a connected dominating set $D$ of $G$. The goal of Dominator is to finish with $|D|$ as small as possible, the goal of Staller is just the opposite. If  both players play optimally, then $|D|$ is unique. In the D-game it is called the {\em connected game domination number} $\cg(G)$ of $G$, while when the S-game is played, the corresponding invariant is denoted by $\cg'(G)$. The connected domination game is thus defined as the standard domination game~\cite{bresar-2010} with the additional requirement that the players maintain connectedness of the subgraph induced by the selected vertices at all times. To shorten the presentation we will abbreviate the term ``connected domination game" to \emph{c}-\emph{game}.

The \emph{connected domination game with Chooser}~\cite{borowiecki-2019} has similar rules as the normal game, except that there is another player, Chooser, who can make zero, one, or more moves after any move of Dominator or Staller. The conditions for his move to be legal are the same as for Dominator and Staller. Chooser has no specific goal, he can help Dominator or Staller or none. We recall the Chooser Lemma from~\cite{borowiecki-2019} to be used later on.

\begin{lemma}[Chooser Lemma]
\label{lema:chooser}
Consider the connected domination game with Chooser on a graph $G$. Suppose that in the game Chooser picks $k$ vertices, and that both Dominator and Staller play optimally. Then at the end of the game the number of played vertices is at most $\cg(G) + k$ and at least $\cg(G) - k$.
\end{lemma}

We now introduce the total connected domination game. First, we recall that the \emph{total domination game} is defined analogously as the domination game, except that whenever a player selects a new vertex in the course of the game, the selected vertex must totally dominate at least one vertex that was not totally dominated by vertices previously selected by the players~\cite{he-kl-ra-2015}. The \emph{game total domination number} and the \emph{Staller-start game total domination number} are denoted by $\ggt(G)$ and $\ggts(G)$, respectively. The \emph{total connected domination game} is just as the total domination game with the additional requirement that at each stage of the game the selected vertices induce a connected subgraph of $G$. The \emph{game total connected domination number} and the \emph{Staller-start game total connected domination number} are denoted by $\tcg(G)$ and $\tcg'(G)$, respectively. For simplicity, we abbreviate the term ``total connected domination game" to \emph{tc}-\emph{game}.

The Chooser Lemma holds also for the total connected domination game. Its proof proceeds along the same lines as the proof of the Chooser Lemma in~\cite{borowiecki-2019}, hence we do not repeat it here.

If $\cc(G) = 1$, then $G$ contains a universal vertex and hence $\tcg(G) = 2$. For all the other cases we have the following bounds that can be proved along the same lines as~\cite[Theorem 1]{borowiecki-2019}, see also~\cite[Theorem 2.1]{bujtas-2019} for more detailed arguments.

\begin{proposition}
\label{prop:bounds-tcc}
If $\cc(G)\ge 2$, then $\cc(G) \le \tcg(G) \le 2\cc(G) - 1$.
\end{proposition}

No matter which game is played, we adopt the convention that the consecutive moves of Dominator are denoted by $d_1, d_2, \ldots$, and the consecutive moves of Staller by $s_1, s_2, \ldots$.

%%%%%%%%%%%%%%%%%%%%%%%%%%%%%%%%%%%%%%%%%%%%%%%%%%%%%
%%%%%%%%%%%%%%%%%%%%%%%%%%%%%%%%%%%%%%%%%%%%%%%%%%%%%
\section{Relating $\tcg(G)$ to $\cg(G)$}
\label{sec:relation}
%%%%%%%%%%%%%%%%%%%%%%%%%%%%%%%%%%%%%%%%%%%%%%%%%%%%%
%%%%%%%%%%%%%%%%%%%%%%%%%%%%%%%%%%%%%%%%%%%%%%%%%%%%%

The main result of this section reads as follows.

\begin{theorem}
\label{thm:only-3-classes}
If $G$ is a connected graph, then $\cg(G) \le \tcg(G) \le \cg(G) + 2$.
\end{theorem}
\proof
To prove the lower bound, let a c-game be played on a graph $G$, call it the {\em R-game} (where R stands for ``real"). In parallel, Dominator imagines that also a tc-game is played on $G$, call it the {\em I-game} (where I stands for ``imagined"). In the R-game Staller plays optimally (and Dominator maybe not), while in the I-game Dominator plays optimally (and Staller maybe not). At the beginning, Dominator selects an optimal first move in the I-game and copies this move to the real game. (The first move of Dominator in the R-game might not be optimal.) After Staller replies with her optimal move in the R-game, Dominator copies her move into the I-game. Note that this move is legal in the I-game. Afterwards Dominator replies optimally in the I-game, and copies his move into the R-game. The two games proceed along the same lines until the games are finished. By the strategy of Dominator, the sequences of moves played in the R-game and in the I-game are the same. Consequently, the number of moves, say $s$, is the same in both games. Since Staller played optimally in the R-game (which is a c-game) we have $\cg(G) \le s$, and since Dominator was playing optimally in the I-game (which is a tc-game) we have $\tcg(G) \ge s$. Hence, $\cg(G) \le s \le \tcg(G)$, proving the left inequality.

In order to prove the right inequality, let again a c-game, called R-game, be played on a graph $G$. Now Staller imagines that a parallel tc-game, called I-game, is played on $G$. In this set-up, in the R-game Dominator plays optimally (and Staller maybe not), while in the I-game Staller plays optimally (and Dominator maybe not). The R-game starts with an optimal move $d_1$ of Dominator, and Staller copies this move into the I-game. Then Staller optimally replies in the I-game with the move $s_1$.  We now distinguish two cases.

\medskip
\emph{Case 1: $s_1$ is a legal move in the R-game.} In this case Staller copies it into the R-game and the game continues along the same lines. As above we infer that the sequences of moves played in the R-game and in the I-game are the same, let $s$ be the number of them. Since Dominator played optimally in the R-game we have $\cg(G) \ge s$, and since Staller was playing optimally in the I-game we have $\tcg(G) \le s$. So, $\tcg(G) \le s \le \cg(G)$.

\medskip
\emph{Case 2: $s_1$ is not a legal move in the R-game.} This situation has happened because $N_G[s_1] \subseteq N_G[d_1]$. In the R-game now, after the move $d_1$, Chooser plays an arbitrary legal move $x$. Staller then imagines in the I-game that $x$ is the second move of Dominator. Note that $x$ is a legal move (of Dominator) in the I-game because in the R-game Chooser dominated at least one vertex not adjacent to $d_1$. Till this moment the sequence of moves played in the R-game is $d_1$, $x$, and the sequence of moves played in the I-game is $d_1$, $s_1$, $x$. In this way the set of vertices dominated in both games is the same. Moreover, in both games it is Staller's turn. Afterwards, both games continue by copying each Staller's move from the I-game to the R-game and by copying each Dominator's move from the R-game to the I-game. Hence in the rest of the games the sequence of moves is the same in both of them. Setting $s$ to be the number of moves played at the end of the R-game, the number of moves played in the I-game is therefore $s+1$. Since the R-game is a c-game in which Dominator played optimally and Chooser  played one move, the Chooser Lemma implies that $\cg(G) \ge s - 1$. On the other hand, since Staller played optimally in the I-game, $\tcg(G) \le s + 1$. Therefore, $\tcg(G) \le s + 1 \le \cg(G) + 2$.~\qed

In view of Theorem~\ref{thm:only-3-classes}, we say that a graph $G$ is
\begin{itemize}
\item {\em Class~$0$}, if $\tcg(G) = \cg(G)$,
\item {\em Class~$1$}, if $\tcg(G) = \cg(G) + 1$, and
\item {\em Class~$2$}, if $\tcg(G) = \cg(G) + 2$.
\end{itemize}

%%%%%%%%%%%%%%%%%%%%%%%%%%%%%%%%%%%%%%%%%%%%%%%%%%%%%
%%%%%%%%%%%%%%%%%%%%%%%%%%%%%%%%%%%%%%%%%%%%%%%%%%%%%
\section{Families of Class $0$, $1$, and $2$ graphs}
\label{sec:classes}
%%%%%%%%%%%%%%%%%%%%%%%%%%%%%%%%%%%%%%%%%%%%%%%%%%%%%
%%%%%%%%%%%%%%%%%%%%%%%%%%%%%%%%%%%%%%%%%%%%%%%%%%%%%

To describe a large family of Class $0$ graphs, we say that the neighborhoods of a graph $G$ are {\em non-inclusive} if for every pair $u$ and $v$ of distinct vertices of $G$ we have $N[u]\not\subseteq N[v]$. Note that the latter condition is trivially fulfilled if $uv\notin E(G)$, hence an equivalent way to say that the neighborhoods of $G$ are non-inclusive is that for every edge $uv\in E(G)$ we have $N(u) \setminus \{v\} \not\subseteq N(v)$.

\begin{proposition}
\label{prp:not-included-neighborhood}
If the neighborhoods of a connected graph $G$ are non-inclusive, then $G$ is Class~$0$.
\end{proposition}
\proof
This result can be deduced from the proof of Theorem~\ref{thm:only-3-classes} as follows. When proving that $\tcg(G)\le \cg(G) + 2$, the condition that the neighborhoods of $G$ are non-inclusive implies that the move $s_1$ of Staller in the I-game is a legal move in the R-game, hence only Case 1 applies. Therefore, $\tcg(G) \le \cg(G)$ holds and because $\cg(G)\le \tcg(G)$ by the first part of the proof of Theorem~\ref{thm:only-3-classes}, we conclude that $\tcg(G) = \cg(G)$. \qed

\vskip -0.25 cm
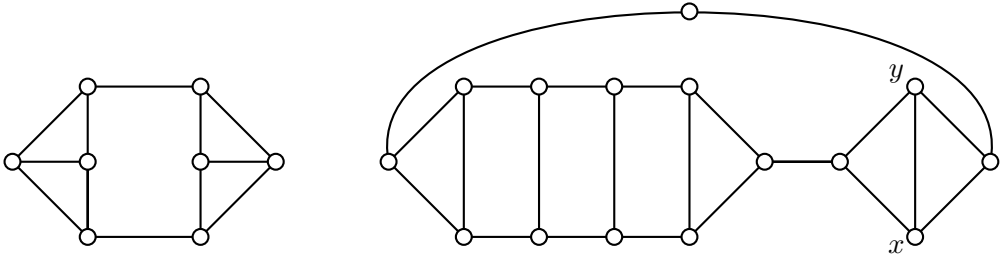
\begin{figure}[t!]
	\begin{center}
		\begin{tikzpicture}[scale=0.5,style=thick]
		\def\vr{6pt}
		%% vertices defined %%
		\begin{scope}
	 \path (0,2) coordinate (b1);
	 \path (2,0) coordinate (a1);
	 \path (2,2) coordinate (a2);
	 \path (2,4) coordinate (a3);
	 \path (5,0) coordinate (a4);
	 \path (5,2) coordinate (a5);
	 \path (5,4) coordinate (a6);
	 \path (7,2) coordinate (b2);

 \draw (b1) -- (a1) -- (a2) -- (a3) -- (a6) -- (a5) --(a4) -- (a1) -- (a2) -- (b1);
	  \draw (b1) -- (a3);
	   \draw (b2) -- (a4);
	    \draw (b2) -- (a5);
	     \draw (b2) -- (a6);
	
	 \draw (b1)  [fill=white] circle (\vr);
	 \draw (a1)  [fill=white] circle (\vr);
	 \draw (a2)  [fill=white] circle (\vr);
	 \draw (a3)  [fill=white] circle (\vr);
	 \draw (a4)  [fill=white] circle (\vr);
	 \draw (a5)  [fill=white] circle (\vr);
	 \draw (a6)  [fill=white] circle (\vr);
	 \draw (b2)  [fill=white] circle (\vr);
		
		\end{scope}
			
		\begin{scope}[xshift = 10cm]
		\path (0,2) coordinate (z1);
		\path (2,0) coordinate (x1);
		\path (2,4) coordinate (u1);
		\path (4,0) coordinate (x2);
		\path (4,4) coordinate (u2);
		\path (6,0) coordinate (x3);
		\path (6,4) coordinate (u3);
		\path (8,0) coordinate (x4);
		\path (8,4) coordinate (u4);
		\path (10,2) coordinate (z2);
		\path (12,2) coordinate (y1);
		\path (14,0) coordinate (y2);
		\path (14,4) coordinate (v2);
		\path (16,2) coordinate (y3);
		\path (8,6) coordinate (z3);
				%% edges %%
		\draw (z1) -- (x1) -- (x2) -- (x3) -- (x4) -- (z2) --(y1) -- (y2) -- (y3) -- (v2) -- (y1) -- (z2) -- (u4) -- (u3) -- (u2) -- (u1) -- (z1);
		\draw (x1) -- (u1);
		\draw (x2) -- (u2);
		\draw (x3) -- (u3);
		\draw (x4) -- (u4);
		\draw (y2) -- (v2);
		 \draw (z1) .. controls (-1,7.3) and (17,7.3) .. (y3);
		%% vertices %%
		\draw (x1)  [fill=white] circle (\vr);
		\draw (y1)  [fill=white] circle (\vr);
		\draw (y2)  [fill=white] circle (\vr);
		\draw (z1)  [fill=white] circle (\vr);
		\draw (x2)  [fill=white] circle (\vr);
		\draw (x3)  [fill=white] circle (\vr);
		\draw (x4)  [fill=white] circle (\vr);
		\draw (y3)  [fill=white] circle (\vr);
		\draw (z3)  [fill=white] circle (\vr);
		\draw (z2)  [fill=white] circle (\vr);
		\draw (u1)  [fill=white] circle (\vr);
		\draw (u2)  [fill=white] circle (\vr);
		\draw (u3)  [fill=white] circle (\vr);
		\draw (u4)  [fill=white] circle (\vr);
		\draw (v2)  [fill=white] circle (\vr);
		%% others %%
		\draw [below] (13.5,0.2) node {{\small $x$}};
		\draw [above] (13.5,3.8) node {{\small $y$}};
				\end{scope}
				\end{tikzpicture}
	\end{center}
\vskip -0.5 cm
	\caption{Two graphs, $F_8$ and $D_{15}$, whose neighborhoods are not non-inclusive, but are Class~$0$}
	\label{fig:class-0}
\end{figure}

The two graphs in Fig.~\ref{fig:class-0} show that the converse of Proposition~\ref{prp:not-included-neighborhood} is not true. $F_8$ is a cubic graph that contains twin vertices, that is, vertices with the same closed neighborhoods. However, as Dominator may choose an optimal start-vertex from the $6$-cycle that does not have a twin in $F_8$, we have $\gamma_{\rm cg}(F_8)= \gamma_{\rm tcg}(F_8)=4$. Therefore, $F_8$ is Class~$0$. More generally, let $F_{4k}$ be the analogous construction on $4k$ vertices. That is, $F_{4k}$ is obtained from a cycle $C_{3k}$ by adding a twin vertex to every third vertex of the cycle. It is clear that $F_{4k}$ is Class~$0$ for each integer $k \ge 2$.

The other graph $D_{15}$ illustrated in Fig.~\ref{fig:class-0} is also Class~$0$ and its neighborhoods are not non-inclusive since $x$ and $y$ are twins. For every $v \in V(D_{15})$ define the invariants $c(v)$ and $t(v)$ as the number of played vertices in a c-game and tc-game, respectively, where Dominator starts the game by playing $d_1=v$ and, after this move the players follow optimal strategies. Then, $\gamma_{\rm cg}(D_{15})= \min_{v\in V(D_{15})} c(v)$ and $\gamma_{\rm tcg}(D_{15})= \min_{v\in V(D_{15})} t(v)$. It can be checked that $\gamma_{\rm cg}(D_{15}) = \gamma_{\rm tcg}(D_{15}) = 9$, hence $D_{15}$ is Class~$0$. Clearly, if $v \in V(D_{15}) \setminus \{x,y\}$, then $c(v)=t(v)$. The interesting fact here is that it can be checked that  $c(x)=t(x)= 10$ and $c(y)=t(y)=10$ also holds, that is, Staller cannot gain any advantage from the fact that she can play a twin in her first move.

By the structure of the Cartesian and the direct product of graphs,  Proposition~\ref{prp:not-included-neighborhood} yields the following two consequences. With respect to the second one we recall that the direct product $G\times H$ is connected if and only if both $G$ and $H$ are connected and at least one of them contains an odd cycle~\cite{hik-2011}. Recall that a non-trivial graph has at least two vertices.

\begin{corollary}
\label{cor:Cartesian-direct}
If $G$ and $H$ are non-trivial graphs, then the following holds.
\begin{enumerate}
\item If both $G$ and $H$ are connected, then $G\cp H$ is Class~$0$.
\item If $\delta(G)\ge 2$ and $G\times H$ is connected, then $G\times H$ is Class~$0$.
\end{enumerate}
\end{corollary}
\proof
(a) Let $(g,h)(g',h)\in E(G\cp H)$. Since $H$ is non-trivial, there exists an edge $hh'\in E(H)$. Then $(g,h')\in N_{G\cp H}((g,h))\setminus N_{G\cp H}((g',h))$. A parallel conclusion can be obtained for an edge $(g,h)(g,h')\in E(G\cp H)$. Hence the neighborhoods of $G\cp H$ are non-inclusive and Proposition~\ref{prp:not-included-neighborhood} applies.

(b) Consider an edge $(g,h)(g',h')\in E(G\times H)$. Since $\delta(G) \ge 2$, there exists an edge $gg''\in E(G)$, where $g''\ne g'$. Now $(g'',h')\in N_{G\times H}((g,h))\setminus N_{G\times H}((g',h'))$. Hence Proposition~\ref{prp:not-included-neighborhood} applies again.
\qed

\medskip 
There exist also direct product graphs that are Class~$0$ but not covered by Corollary~\ref{cor:Cartesian-direct}(b). Let $H$ be the graph obtained from $K_{1,3}$ by adding an arbitrary edge to it ($H$ is known as the paw graph). Then, $\cg(H \times K_2) = \tcg(H \times K_2) = 5$ and, of course, $\delta(H) = \delta(K_2) = 1$. Moreover, $C_{2k+1}\times K_2 = C_{4k+2}$ is Class~0 for every $k\ge 1$. But not all direct products are Class~0. For instance, it can be easily checked that
$$\cg((C_{2k+1}\odot K_1) \times K_2) = 4k+2$$ and
$$\tcg((C_{2k+1}\odot K_1) \times K_2) = 4k+3\,.$$

We next show that generalized coronas over connected graphs are Class~$1$.

\begin{proposition}
\label{prp:generalized-coronas}
If $G$ is a connected graph and ${\cal H} = \{H_1,\ldots, H_{n(G)}\}$, then $G\odot {\cal H}$ is Class~$1$.
\end{proposition}
\proof
Observe first that $\cc(G\odot H) = n(G)$. Therefore, $\tcg(G \odot H) \ge \cg(G \odot H) \ge n(G)$. If in the c-game Dominator selects as his first move a vertex of $G$, then Staller must reply with an adjacent vertex of $G$. Proceeding in this way Dominator can ensure that exactly all the vertices of $G$ will be played by the end of the game, hence $\cg(G\odot {\cal H}) \le n(G)$ and so $\cg(G\odot {\cal H}) = n(G)$. In the tc-game, Dominator's optimal strategy is again to play a vertex $v_i$ of $G$. If Staller replies with a vertex of $G$, then, as above, only the vertices of $G$ will be played. Hence an optimal reply of Staller is to play a vertex from $H_i$. After that Dominator can ensure that only the remaining vertices of $G$ will be played, so that $\cg(G\odot {\cal H}) = n(G) + 1$. We conclude that $G\odot {\cal H}$ is Class~$1$.
\qed

\begin{corollary}
\label{cor:trees}
A tree $T$ with $n(T)\ge 3$ is Class~$1$ if and only if $T = T'\odot {\cal H}$, where $T'$ is a tree and ${\cal H}$ is a family of edge-less graphs. Otherwise, $T$ is Class~$0$.
\end{corollary}
\proof
As observed in~\cite{borowiecki-2019}, $\cg(T) = n(T) - \ell(T)$, where $\ell(T)$ is the number of leaves of $T$. If $T$ contains a vertex $x$ of degree at least $2$ with no leaf attached to it, then Dominator plays $x$ first and in this way guarantees that exactly the non-leaves will be played. Hence, $T$ is Class~$0$ in this case. Otherwise every non-leaf has at least one leaf attached. In this case, $T$ can be represented as $T'\odot {\cal H}$ for some tree $T'$ and a family ${\cal H}$ of edge-less graphs. By Proposition~\ref{prp:generalized-coronas}, $T$ is Class~$1$ in this case.
\qed

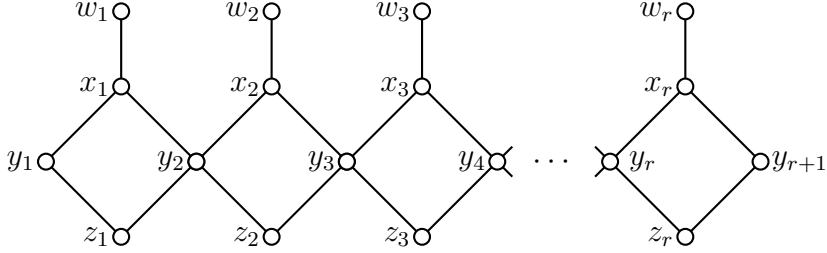
\begin{figure}[ht!]
\begin{center}
\begin{tikzpicture}[scale=1.0,style=thick]
\def\vr{3pt}
%% vertices defined %%
\begin{scope}
\path (0,1) coordinate (y1);
\path (1,0) coordinate (z1);
\path (1,2) coordinate (x1);
\path (1,3) coordinate (w1);
\path (2,1) coordinate (y2);
%% edges %%
\draw (w1) -- (x1) -- (y1) -- (z1) -- (y2) -- (x1);
%% vertices %%
\draw (x1)  [fill=white] circle (\vr);
\draw (y1)  [fill=white] circle (\vr);
\draw (y2)  [fill=white] circle (\vr);
\draw (z1)  [fill=white] circle (\vr);
\draw (w1)  [fill=white] circle (\vr);
%% others %%
\draw [left] (y1) node {$y_1$};
\draw [left] (z1) node {$z_1$};
\draw [left] (y2) node {$y_2$};
\draw [left] (x1) node {$x_1$};
\draw [left] (w1) node {$w_1$};
\end{scope}
%%%%%%%%%%%%%%%%%%%%%%%%%%%%%%%%%%%%
\begin{scope}[xshift = 2cm]
\path (0,1) coordinate (y1);
\path (1,0) coordinate (z1);
\path (1,2) coordinate (x1);
\path (1,3) coordinate (w1);
\path (2,1) coordinate (y2);
%% edges %%
\draw (w1) -- (x1) -- (y1) -- (z1) -- (y2) -- (x1);
%% vertices %%
\draw (x1)  [fill=white] circle (\vr);
\draw (y1)  [fill=white] circle (\vr);
\draw (y2)  [fill=white] circle (\vr);
\draw (z1)  [fill=white] circle (\vr);
\draw (w1)  [fill=white] circle (\vr);
%% others %%
\draw [left] (z1) node {$z_2$};
\draw [left] (y2) node {$y_3$};
\draw [left] (x1) node {$x_2$};
\draw [left] (w1) node {$w_2$};
\end{scope}
%%%%%%%%%%%%%%%%%%%%%%%%%%%%%%%%%%%%
\begin{scope}[xshift = 4cm]
\path (0,1) coordinate (y1);
\path (1,0) coordinate (z1);
\path (1,2) coordinate (x1);
\path (1,3) coordinate (w1);
\path (2,1) coordinate (y2);
%% edges %%
\draw (w1) -- (x1) -- (y1) -- (z1) -- (y2) -- (x1);
\draw (2.2,1.2) -- (2,1) -- (2.2,0.8);
%% vertices %%
\draw (x1)  [fill=white] circle (\vr);
\draw (y1)  [fill=white] circle (\vr);
\draw (y2)  [fill=white] circle (\vr);
\draw (z1)  [fill=white] circle (\vr);
\draw (w1)  [fill=white] circle (\vr);
%% others %%
\draw [left] (z1) node {$z_3$};
\draw [left] (y2) node {$y_4$};
\draw [left] (x1) node {$x_3$};
\draw [left] (w1) node {$w_3$};
\end{scope}
%%%%%%%%%%%%%%%%%%%%%%%%%%%%%%%%%%%%
\begin{scope}[xshift = 7.5cm]
\path (0,1) coordinate (y1);
\path (1,0) coordinate (z1);
\path (1,2) coordinate (x1);
\path (1,3) coordinate (w1);
\path (2,1) coordinate (y2);
%% edges %%
\draw (w1) -- (x1) -- (y1) -- (z1) -- (y2) -- (x1);
\draw (-0.2,1.2) -- (0,1) -- (-0.2,0.8);
%% vertices %%
\draw (x1)  [fill=white] circle (\vr);
\draw (y1)  [fill=white] circle (\vr);
\draw (y2)  [fill=white] circle (\vr);
\draw (z1)  [fill=white] circle (\vr);
\draw (w1)  [fill=white] circle (\vr);
%% others %%
\draw [left] (z1) node {$z_r$};
\draw [right] (y2) node {$y_{r+1}$};
\draw [right] (y1)++(0.1,0) node {$y_{r}$};
\draw [left] (x1) node {$x_r$};
\draw [left] (w1) node {$w_r$};
\end{scope}
\draw (6.75,1) node {$\cdots$};

\end{tikzpicture}
\end{center}
\caption{The graph $G_r$}
\label{fig:cacti-class-2}
\end{figure}

Let $G_r$, $r\ge 3$, be a class of graphs whose formal definition should be clear from Fig.~\ref{fig:cacti-class-2}.

\begin{proposition}
\label{prp:cacti-class-2}
If $r\ge 3$, then $G_r$ is Class~$2$.
\end{proposition}
\proof
Let $r\ge 3$ and let the vertices of $G_r$ be labeled as shown in Fig.~\ref{fig:cacti-class-2}. Let $X= \{x_1,x_2, \ldots, x_r\}$ and let $Y = \{y_2,y_3, \ldots, y_r\}$. Every connected dominating set of $G_r$ contains the set $X \cup Y$. On the other hand, $X \cup Y$ is itself a connected dominating set of $G_r$. Consequently, $\cc(G_r) = |X \cup Y| = 2r-1$ (and the set $X \cup Y$ is the unique minimum connected dominating set of $G_r$). Hence, $\cg(G_r) \ge \cc(G_r) = 2r-1$.

We first consider the c-game played on $G_r$. Set $d_1 = x_2$. Then, $s_1\in \{y_2, y_3\}$. The strategy of Dominator in the rest of the game is the following. Suppose by induction that $y_i$, $2\le i\le r$, is the currently last vertex played by Staller. Then Dominator plays the vertex from $\{x_{i-1}, x_{i}\}$ that has not yet been played. Maintaining this strategy, Dominator ensures that Staller is forced always to play a vertex from $Y$. Moreover, by the end of the game, Dominator will play all vertices from $X$, so that $\cg(G_r) \le 2r-1$ and hence $\cg(G_r) = 2r-1$.

We consider next the tc-game on $G_r$, and show that $\tcg(G) \ge 2r + 1 = \cg(G) + 2$, implying that $G_r$ is Class~$2$. As observed earlier, every connected dominating set of $G_r$ contains the set $X \cup Y$. Hence, it suffices for us to show that Staller has a strategy that forces at least two vertices of $G$ to be played that do not belong to the set $X \cup Y$. In this case, we say that Staller achieves her goal. Let $V_i = \{w_i,x_i,y_i,y_{i+1},z_i\}$ where $i \in [r]$.

%If Dominator plays as his first move $d_1 = y_1$, then Staller responds with $s_1 = z_1$, while if $d_1 = z_1$, then Staller responds with $s_1 = y_1$. In both cases, two moves not in $X \cup Y$ are played, as desired. Hence, we may assume that $d_1 \notin \{y_1,z_1\}$. By symmetry, we may assume that $d_1 \notin \{y_{r+1},z_r\}$. Let $V_i = \{w_i,x_i,y_i,y_{i+1},z_i\}$ where $i \in [r]$.

Suppose that $d_1 \in V_1$. If $d_1 = y_1$, then Staller plays $s_1 = z_1$, while if $d_1 = z_1$, then Staller plays $s_1 = y_1$. In both cases, two moves not in $X \cup Y$ are played, and Staller achieves her goal. If $d_1 = w_1$, then $s_1 = x_1$, and either $d_2 = y_1$ or $d_2 = y_2$. If $d_2 = y_1$, then Staller immediately achieves her goal, while if $d_2 = y_2$, then Staller plays $s_2 = z_2$ to achieve her goal. If $d_1 = x_1$, then Staller plays $s_1 = w_1$, and either $d_2 = y_1$ or $d_2 = y_2$, and as in the previous case Staller can achieve her goal. Hence, Staller achieves her goal that two vertices are played from outside the set $X \cup Y$ in all cases except possibly in the case when $d_1 = y_2$. Suppose, therefore, that $d_1 = y_2$. In this case, Staller plays $s_1 = z_2$. If $d_2 = z_1$, then Staller immediately achieves her goal. If $d_2 \in \{x_2,y_3\}$, then Staller plays $s_2 = z_1$, and achieves her goal. If $d_2 = x_1$, then Staller plays $s_2 = x_2$, which forces $d_3 = y_3$, and enables Staller to play $s_3 = z_3$, and again she achieves her goal. Hence, if $d_1 \in V_1$, then Staller achieved her goal. We may therefore assume that $d_1 \notin V_1$ and, by symmetry, that $d_1 \notin V_r$. Hence, $d_1 \in V_i$ for some $i \in [r-1] \setminus \{1\}$.

Suppose that $d_1 = w_i$, which forces $s_1 = x_i$, and either $d_2 = y_i$ or $d_2 = y_{i+1}$. If $d_2 = y_i$, then Staller plays $s_2 = y_{i+1}$, while if $d_2 = y_{i+1}$, then Staller plays $s_2 = y_{i}$. Thus, $d_3 \in \{x_{i-1},z_{i-1},x_{i+1},z_{i+1}\}$. If $d_3 \in \{x_{i-1},z_{i-1}\}$, then Staller plays $s_3 = z_{i+1}$, while if $d_3 \in \{x_{i+1},z_{i+1}\}$, then Staller plays $s_3 = z_{i-1}$. In both cases, after her third move Staller already achieve her goal that two vertices are played from outside the set $X \cup Y$.

Suppose that $d_1 = x_i$. In this case, Staller plays $s_1 = w_i$, which forces either $d_2 = y_i$ or $d_2 = y_{i+1}$. Proceeding exactly as in the previous case, if $d_2 = y_i$, then Staller plays $s_2 = y_{i+1}$, while if $d_2 = y_{i+1}$, then Staller plays $s_2 = y_{i}$, thereby achieving her goal as before.

Suppose that $d_1 = y_i$. In this case, Staller plays $s_1 = z_i$. We note that $d_2 \in \{x_{i-1},z_{i-1},x_{i},y_{i+1}\}$. If $d_2 = z_{i-1}$, then already two vertices are played from outside the set $X \cup Y$, and Staller immediately achieves her goal. If $d_2 \in \{x_{i},y_{i+1}\}$, then Staller plays $s_2 = z_{i-1}$ and achieves her goal. Hence, we may assume that $d_2 = x_{i-1}$. In this case, Staller plays $s_2 = x_i$, forcing either $d_3 = y_{i-1}$ or $d_3 = y_{i+1}$. If $d_3 = y_{i-1}$ and $i = 2$, then Staller immediately achieves her goal. If $d_3 = y_{i-1}$ and $i \ge 3$, then Staller plays $s_3 = z_{i-2}$, while if $d_3 = y_{i+1}$, then Staller plays $s_3 = z_{i+1}$, and in both cases she achieves her goal. Analogously, if $d_1 = y_{i+1}$, then Staller achieves her goal.

Suppose finally that $d_1 = z_i$. In this case, Staller plays $s_1 = y_i$. Thus, $d_2 \in \{x_{i-1},z_{i-1},x_{i},y_{i+1}\}$. If $d_2 = z_{i-1}$, then Staller immediately achieves her goal. If $d_2 \in \{x_{i},y_{i+1}\}$, then Staller plays $s_2 = z_{i-1}$ and achieves her goal. Hence, we may assume that $d_2 = x_{i-1}$. In this case, Staller plays $s_2 = x_i$, forcing either $d_3 = y_{i-1}$ or $d_3 = y_{i+1}$. Staller now proceeds therefore exactly as in the previous case to achieve her goal. \qed

%%%%%%%%%%%%%%%%%%%%%%%%%%%%%%%%%%%%%%%%%%%%%%%%%%%%%
%%%%%%%%%%%%%%%%%%%%%%%%%%%%%%%%%%%%%%%%%%%%%%%%%%%%%
\section{Concluding remarks}
\label{sec:problems}
%%%%%%%%%%%%%%%%%%%%%%%%%%%%%%%%%%%%%%%%%%%%%%%%%%%%%
%%%%%%%%%%%%%%%%%%%%%%%%%%%%%%%%%%%%%%%%%%%%%%%%%%%%%

In Proposition~\ref{prop:bounds-tcc} the bounds are sharp. The lower bound is attained by trees which are Class~$0$, while the upper bound is attained by classes of Cartesian products $X$ from~\cite{borowiecki-2019, bujtas-2019} for which $\cg(X) = 2\cc(X) - 1$. Hence the upper bound in Proposition~\ref{prop:bounds-tcc} is  attained because Cartesian products are Class~$0$ graphs.

\begin{problem}
Determine whether in Proposition~\ref{prop:bounds-tcc} all possible values of $\tcg(G)$ are realizable.
\end{problem}

With respect to Corollary~\ref{cor:Cartesian-direct}(b) we pose:

\begin{problem}
Classify direct product graphs into Classes $0$, $1$, and $2$.
\end{problem}

Analyzing the graphs $G_n$ from~\cite{irsic-2019+} (see also Fig.~1 there), it can be demonstrated that for every natural number $k$  there exist graphs $G$ such that $\tcg'(G) - \tcg(G) \ge k$. In this paper we do not however further investigate the S-game, hence the following task remains to be done.

\begin{problem}
Consider the total connected domination S-game. In particular, we suspect that $\tcg'(G) = \cg'(G)$ holds whenever $G$ is not a complete graph.
\end{problem}

We conclude with the following problem that also seems to be interesting.

\begin{problem}
Classify cactus graphs into Classes $0$, $1$, and $2$.
\end{problem}

\section*{Acknowledgements}

We acknowledge the financial support from the Slovenian Research Agency (research core funding No.\ P1-0297 and 
projects J1-9109, J1-1693, N1-0095, N1-0108).

\end{document}